\documentstyle[11pt]{article}
\begin{document}
\title{ CONSTRUCTION OF SEMI - FREE $S^3$\\
ACTIONS ON HOMOTOPY
SPHERE WITH UNTWISTED FIXED POINT SET}
\vspace{4cm}
\author{Issam. H. Kaddoura.\footnote{e.mail~:~issam.kaddoura@liu.edu.lb.}~\footnote{
Faculty of Science Baghdad University. } \\
Lebanese International University \\
Department of Mathematics
 }
\date{}
\maketitle
\begin{abstract}
William Browder in his paper "Surgery and the theory of differentiable
transformation groups " developed surgery techniques to study semi-free actions of
$S^1$ on homotopy spheres, under the additional assumption that the fixed point set is a
homotopy sphere. He used this surgery to show how to construct such actions. In this
paper, I discussed a similar theory for semi-free actions of $S^3$ on homotopy spheres.\\
An open problem is raised at the end of the paper.
\end{abstract}
{\it Keywords:} Free and semi free actions, quaterionic projective space, homotopy and
homology spheres, principal fibration, tangent and normal bundles,
exact homology and cohomology sequences, lefshetz duality,
Poincare duality, Hurwicz homomorphism, Hopf fibration, h~-~cobordism,
equivariant diffeomorphism, VanKampen�s theorem, Mayer~-~Vietoris sequence,
untwisted fixed point set, K\"{u}nneth formula, surgery techniques.
\section{Introduction }
 Through this paper, $R^n$ denotes the Euclidean $n-$space, $S^n$ denotes the unit
$n-$sphere in $R^{n+1}$ and $QP(n)$ the quaterionic projective space, all having the usual
differentiable structures. By a homotopy $n-$sphere abbreviated by $\sum^n$, we mean a
closed differentiable $n-$manifold having the homotopy type of $S^n$, and by a homotopy
quaterionic projective n-space,abbreviated by $HQP(n)$, we mean a closed
differentiable $4n-$manifold having the homotopy type of $QP(n)\cdot \large{\pi}_n(M)$ denotes the
$n^{th}-$homotopy group of $M$, $H_i(M,G)$ and $H^i(M,G)$ denote the homology and cohomology
of a space $M$, with coefficients in the group $G$ and assumed to be satisfying the Eeilenberg$-$Steenrod axioms, See
[\cite{14},p.6]. If $G = Z$, we write $H^i ( M)$ and $H_i (M)$
for $H^i ( M,Z )$ and $H_i (M,Z)$ respectively. It is well known that $S^1$ and $S^3$ are the only compact connected Lie groups
which have free differentiable actions on homotopy spheres \cite{4,5}. It follows from
Gleason's lemma \cite{1} that such an action is always a principal fibration which is
homotopically equivalent to the classical Hopf fibration.\\In fact, there are always
infinitely many differentiably distinct free actions of $S^3$ on $\sum^{4n+3}$ for $n \geq 2$ \cite{16}.

\section{Construction of  semi$-$free actions of $S^3$:}
An action $(G, M, \phi )$ is called semi$-$free, if it
is free outside the fixed point set, that is $\phi: G\times M \rightarrow M$, and $F$ is the fixed point set of the
action, $\phi$ is semi$-$free action if $\phi(g,x)= x$, for some $x \in M~-~F$ then $g = e$ the
identity of $G$. Notice that, there are two types of orbits, fixed points, and $G$.\\\\
{\bf  Lemma 2.1 } Let $\phi : S^i \times M \rightarrow M~~~( i= 1,3),$ be a semi free differentiable action. let $F^k$ denote the
union of the $k-$dimensional components of the set of all fixed points of $\phi$. Then the
normal bundle of an imbedding $F^k \subset M$ has naturally a complex strucure for $i=1$
and a quaterionic structure for $i=3$ and that the induced $S^i-$action on the normal
bundle is a scalar multiplication. For proof see \cite{3}.\\\\
 It follows from this lemma $2.1$, that the
codimension of each component of $F$ in $M$ is even for $i=1$ and is divisible by $4$ for $i=3$.
 We shall study the situation where $(S^3,\sum^m , \phi )$ is a semi$-$free differentiable
action on a homotopy sphere $\sum^m$ and the fixed point set is a homotopy sphere $\sum^r$.
Let $(S^3, \sum^m, \phi)$ be a semi$-$free action with fixed point set  $\sum^r \subset \sum^m$, then $S^3$
acts freely outside $\sum^r$ and $S^3$ acts freely and linearly on the normal space to $\sum^r$ at
each point of $\sum^r$. See [\cite{9},p.58]. By lemma 2.1, the normal bundle of $\sum^r$ has a
quaterionic structure and $m -r= 4k,~~k \geq 1$. Let $\mu$ be the ( quaterionic) bundle
over $\sum^r$ defined by the action. We prove the following :\\\\
{\bf  Theorem 2.1 }  $~$If $(S^3, \sum^m, \phi_1)$ and $(S^3, \sum^m, \phi_2)$ are equivalent, then $F_1$ is
diffeomorphic to $F_2$ and $\mu_1$ is equivalent to $\mu_2$, where $F_1~ \mbox{and}~ F_2$ are the fixed point sets
of $\phi_1~\mbox{and}~\phi_2$ respectively, $\mu_1~\mbox{and}~\mu_2$ are the normal bundles of $F_1~\mbox{and}F_2$
respectively.\\
{\bf Proof.} $~$ Since $(S^3,\sum^m, \phi_1)$ and $(S^3,\sum^m, \phi_2)$ are equivalent then
there exists an equivariant diffeomorphism f such that the following diagram
commutes :
\begin{eqnarray*}
S^3\times \mbox{$\sum^m$} \buildrel {\phi_1} \over\longrightarrow \mbox{$\sum^m$}\hspace{3cm}\\
I \times f \downarrow ~~~~~~f \downarrow\hspace{3cm}\\
S^3\times \mbox{$\sum^m$} \buildrel {\phi_2} \over\longrightarrow \mbox{$\sum^m$}\hspace{3cm}
\end{eqnarray*}
Let $x \in F_1$, i.e, $\phi_1(q,x)= x~~\forall ~q\in  S^3~~\mbox{then} ~~~f\circ \phi_1(q,x)= \phi_2 \circ (I\times f)(q,x)~~ \forall~ q\in S^3,$
thus $f(x)= \phi_2(q,f(x))~~ \forall~ q \in S^3~~ \mbox{implies}~~ f(x)\in F_2~~ \mbox{hence}~~ f(F_1) \subseteq	F_2$.\\
Now, assume that $y \in F_2$, i.e., $\phi_2(q,y)=y~~\forall~ q\in S^3$, but $f$ is an equivariant diffeomorphism, then
 $\exists ~~ x \in \sum^m~~ \mbox{such that}~~ y=f(x) ~~\mbox{and}~~ f\circ \phi_1(q,x)= \phi_2 \circ (I\times f)(q,x)~~ \forall~ q\in S^3,$ then
 $f( \phi_1(q,x))= \phi_2 (q,f(x))=f(x),~~\mbox{but}~ f ~\mbox{is}~ 1-1,~\mbox{ we get}~~  \phi_1(q,x)=x ~~\forall~~q\in S^3,$ hence $x\in  F_1$
and $f(x)\in f(F_1)$, therefore $f(F_1) = F_2$. Moreover, the equivalence
$f:\sum^m \rightarrow \sum^{m}$
defines a quaterionic map of the normal bundles $\mu_1$ and $\mu_2$ of $F_1$ and $F_2$
respectively, so they are equivalent.\\\\
Now, $\sum^r$ is a closed submanifold of $\sum^m$ and
invariant under the action of $S^3$ ,therefore there exists a tubular neighbourhood $E$ of
$\sum^r$ which is invariant under the action of $S^3$, so we may consider $\mu : E \rightarrow \sum^r$, the
normal bundle to $\sum^r$ in $\sum^m$ and $S^3$ acts differentiably on $E$. See [\cite{9},p.58].
Let  $~S^{4k-1}~$ be the boundary of a fibre of $E$,\\
we can prove the following :\\\\
{\bf  Lemma 2.2 } $~$ $~S^{4k-1}~\subset \sum^m -\sum^r$ is a homotopy equivalence.\\
{\bf Proof.} $~$ consider the exact cohomology sequence of the pair $(\sum^m, \sum^r):$\\
\[
\rightarrow H^{i-1}(\mbox{$\sum^r$})\rightarrow H^i(\mbox{$\sum^m$}, \mbox{$\sum^r$})
\rightarrow H^i(\mbox{$\sum^m$})\rightarrow H^i( \mbox{$\sum^r$})\rightarrow H^{i+1}(\mbox{$\sum^m$},\mbox{$\sum^r$})\rightarrow H^{i+1}(\mbox{$\sum^m$})\rightarrow
\]
\begin{itemize}
\item[] At $i = m : \rightarrow 0\rightarrow H^m(\sum^m, \sum^r)\rightarrow Z\rightarrow 0, \mbox{then}~ H^m(\sum^m, \sum^r) \cong Z $,
\item[] At $i = r : 0\rightarrow Z \rightarrow H^{r+1}(\sum^m, \sum^r)\rightarrow 0, \mbox{hence}~  H^{r+1}(\sum^m, \sum^r)\cong Z,$
\item[] And for $i \neq m , r+1 : 0 \rightarrow 0\rightarrow H^i(\sum^m, \sum^r)\rightarrow 0$,
\end{itemize}
thus  $H^i(\sum^m, \sum^r)\cong 0$ for $i\neq m,~r+1.$ \mbox{Finally, we obtain}\\
\[
 H^i(\mbox{$\sum^{m}$}, \mbox{$\sum^{r}$})=
\left\{\begin{array}{ll}
Z~~~\mbox{for}~~i=m,~ r+1
\vspace{.2cm}\\
0 ~~~~~~~~~~~~~\mbox{otherwise}
\end{array}
\right.
\]\\
It is clear that $\sum^m- \sum^r$ is simply connected [\cite{66},p.3] and using Lefshetz duality [\cite{6},p.32],\\
$H_i(\sum^m- \sum^r)\cong H^{m-i}(\sum^m, \sum^r)$, we finally deduce that\\
\[
 H_i(\mbox{$\sum^{m}$}- \mbox{$\sum^{r}$})=
\left\{\begin{array}{ll}
Z~~~\mbox{for}~~~~~ i=0 , 4k-1~~~
\vspace{.2cm}\\
0 ~~~~~~~~~~~~~~~\mbox{otherwise}
\end{array}
\right.
\]
hence, the inclusion map $i : S^{4k-1}\rightarrow  \sum^m- \sum^r$ induces\\
$i_\ast : \pi_j(S^{4k-1})\rightarrow  \pi_j(\sum^m- \sum^r)$ an isomorphism $\forall~j$, see[\cite{12},p.283]. Therefore\\
$S^{4k-1}\subset  \sum^m- \sum^r$ is a homotopy equivalence.\\\\
Now, let $N=\sum^m- E_0$ where $E_0$ is the interior of an equivariant tubular
neighbourhood of $\sum^r$ with $\overline{E_0} \subset int(E)$. Then $S^3$ acts freely on $N$ and $S^{4k-1} \subset N$.
Notice that $S^{4k-1}$ is a homotopoy equivalence to $N$, it follows from the exact
homotopy sequence of the fibre maps, using the diagram :
\begin{eqnarray*}
\rightarrow S^3 \rightarrow ~&S^{4k-1}&~\rightarrow S^{4k-1}/{S^3}\rightarrow \hspace{3cm}\\
\downarrow~~~~~~&\downarrow&~~~~~~~\downarrow\hspace{3cm}\\
 \rightarrow S^3 \rightarrow ~&N&~\rightarrow {N}/{S^3}\rightarrow\hspace{3cm}
\end{eqnarray*}
that is $S^{4k-1}/{S^3}\rightarrow N/S^3$ is a homotopy equivalence.\\
Set $N^1 ={N}/{S^3}$ and $S^{4k-1}/{S^3}=QP(k-1)$. Notice that the region between $\partial {N}^1$ and $QP( k-1)\times S^r$ is an
h$-$cobordism, so if $m \geq 6$, then by the h$-$cobordism theorem of Smale \cite{13},
$N^1$ is diffeomorphic to $QP( k-1)\times D^{r+1}$ and $N\rightarrow  N^1$ is equivalent to\\
$h\times I : S^{4k-1}\times D^{r+1}\rightarrow  QP( k-1)\times D^{r+1}$, where $h : S^{4k-1} \rightarrow QP( k- 1)$ is the hopf map.
Hence, we have proved the following theorem:\\\\
{\bf  Theorem 2.2 } $~$ Let $(S^3, \sum^m, \phi)$ be a semi$-$free $S^3$ action on $\sum^m$ with fixed
point set $\sum^r,~ m-r= 4k,~ k\geq1,~m\geq6.$ If $N$ is the complement of an open tubular
neighbourhood of $\sum^r$ in $\sum^m$,then $N$ is equivariantly diffeomorphic to $S^{4k-1}\times D^{r+1}$,
with the standard action on $S^{4k-1}$ and trivial action on $D^{r+1}$.\\\\
Now, we will describe how to construct smooth semi$-$free $S^3$ actions on a
homotopy $m-$sphere $\sum^m$. Let $\sum^r$ be a homotopy $r-$sphere and $\mu$ a (quaterionic)
normal bundle over $\sum^r$ given by $\mu: E(\mu)\rightarrow \sum^r$ where $E(\mu)$ is the total space of
$\mu$ such that $E(\mu)\cong D^{4k}\times \sum^r$, i.e., $E(\mu)$ is the a trivial bundle and suppose that
$h : S^{4k-1}\times \sum^r \rightarrow  S^{4k-1}\times S^r$ is an equivariant diffeomorphism .\\\\
{\bf Theorem 2.3 } $~$ There is a semi free action $(S^3, \sum^m, \phi)$ with a fixed point set
$\sum^r$ and $\sum^m =E(\mu)\cup_{h} (S^{4k-1}\times D^{r+1}) $ where $\cup_h$ means that we identify
$S^{4k-1}\times \sum^r \subset E(\mu)$ with $S^{4k-1}\times S^r \subset S^{4k-1}\times D^{r+1}$ via the diffeomorphism $h$.\\
{\bf Proof.} $~$ Consider the semi$-$free action on the total space of $\mu:E(\mu)\rightarrow \sum^r$ defined by the quaterionic structure and the free $S^3 -$ action on $S^{4k-1}\times D^{r+1}$
defined by the free action on $S^{4k-1}$, i.e., the standard action and
 $h:S^{4k-1}\times \sum^r\rightarrow S^{4k-1}\times S^r$  is an equivariant diffeomorphism, then $M= E(\mu)\cup_h (S^{4k-1}\times D^{r+1})$ has a semi$-$free action of $S^3$ with fixed point set $\sum^r$ and normal
bundle $\mu$, it is enough to show that $M$ is a homotopy sphere.\\
It is clear that $$\pi_1(\partial E(\mu))\cong \pi_1 (S^{4k-1}\times S^r)\cong \pi_1 (S^{4k-1})\oplus( S^r)\cong \pi_1 (\mbox{$\sum^r$})~~\mbox{and}~~ \pi_0(\partial E(\mu))\cong \pi_0 (\mbox{$\sum^r$})\cong 0,$$i.e., $E(\mu)$ and $S^{4k-1}\times D^{r+1}$ are simply connected and $E(\mu)\cap S^{4k-1}\times D^{r+1}$ is simply connected, hence by VanKampen's theorem \cite{10}, $M$ is simply connected.\\
\\
Now, we consider the Mayer -Vietoris sequence for $M$ \cite{2}:
\[\rightarrow H_{s+1}(M)\rightarrow H_s(\partial E(\mu))\rightarrow H_s( E(\mu))\oplus H_s(S^{4k-1}\times D^{r+1})\rightarrow H_s(M)\rightarrow
\]
By the K$\ddot{\mbox{u}}$nneth formula [\cite{8}, p.98], since $H_s(\partial E(\mu)),~H_s(E(\mu)),~H_s(S^{4k-1}\times D^{r+1})$
are torsion free for,  $~0 < s < 4k+ r-1$, we obtain
\[H_s(\partial E(\mu))=H_s(S^{4k-1}\times \mbox{$\sum^{r}$})\cong \oplus_{i=0}^s H_i(S^{4k-1})\otimes H_{s-i}( \mbox{$\sum^{r}$}).
\]
\begin{itemize}
\item[] If $i=0$, then $H_0(S^{4k-1})\otimes H_s(S^r)\cong Z\otimes H_s(\sum^r)$ .
\item[] If $i=s$, then $H_s(S^{4k-1})\otimes H_0(S^r)\cong H_s(S^{4k-1})\otimes Z$ and for $i\neq 0$
\item[] $H_i(S^{4k-1})\otimes H_{n-i}(S^r)\cong 0$ therefore $H_s(\partial E(\mu))\cong Z\otimes H_s(S^r)\oplus  H_s(S^{4k-1})\otimes Z.$
\end{itemize}
Again, we compute
\begin{itemize}
\item[] $H_s(S^{4k-1}\times D^{r+1})=\oplus_{i=0}^s ~H_i(S^{4k-1})\otimes H_{s-i}(D^{r+1})\cong H_s(S^{4k-1})\otimes Z.$
\end{itemize}
\begin{itemize}
\item[] similarly $H_s(E(\mu))=H_s(D^{4k}\times \sum^r)\cong Z\otimes H_s(\sum^r)$.
\item[] therefore $H_s(\partial E(\mu))\cong H_s( E(\mu))\oplus H_s(S^{4k-1}\times D^{r+1})$,
\item[] hence $H_s(M)\cong H_{s+1}(M) \cong 0,~~\forall ~0< s < 4k+ r-1$.
\end{itemize}
For, $s=4k+r-1~:~~H_{4k+r-1}(E(\mu))\cong 0,~~H_{4k+r-1}(S^{4k-1}\times D^{r+1})\cong 0~~~\mbox{and} \\
H_{4k+r-1}(\partial E(\mu)) \cong Z\otimes Z\cong Z.$\\
\\
Substituting in the Mayer -Vietoris sequence for $M : 0 \rightarrow H_{4k+r}(M)\rightarrow Z \rightarrow 0.$
\mbox{~Finally, we obtain}
\[
 H_s(M)\cong
\left\{\begin{array}{ll}
z~~~ \mbox{for}~~ s= 0, ~~4k+r,
\vspace{.2cm}\\
0~~~~~~~~~~~~~otherwise
\end{array}
\right.
\]
hence $M$ is a homology sphere.\\
But $M$ is simply connected, closed without boundary, therefore $M$ is a homotopy
sphere.
\section{Applying surgery to construct semi free $S^3-$ actions:}
In this section, we used surgery techniques as Browder \cite{15} to create a
diffeomorphism of $QP( k-1)\times \sum^r$ with $QP(k-1)\times S^r$, then we apply Theorem 2.3\\\\
{\bf  Theorem 3.1} $~$ Let $\sum^{4n-1}$ be a homotopy sphere which bounds a parallelizable
manifold, $n\geq 1$. Then for each even $k\geq 2$, there is a semi$-$free action of $S^3$ on a
homotopy sphere $\sum^{4(n+k)-1}$ with $\sum^{4n-1}$ as untwisted fixed point set.\\
{\bf Proof.} $~$ Let $\sum^{4n-1}=\partial W^{4n},~~W$ is a parallelizable manifold. We may consider \\
$W_0=W-int(D^{4n})$ as a parallelizable cobordism between $\sum^{4n-1}$ and $S^{4n-1}$
thus we may define a normal map
\[
f:(W_0,~\mbox{$\sum^{4n-1}$} \cup S^{4n-1}) \rightarrow( S^{4n-1}\times I,~S^{4n-1}\times \{0\}\cup S^{4n-1}\times\{1\})
\]
with $f|_{S^{4n-1}}=$ Identity. Since $W$ is a parallelizable manifold [\cite{7},p.514], we may
assume that $W_0$ is $(2n-1)$ connected.\\
Multiplying by $QP(k-1)$, we get $I\times f$:
\[
Qp(k-1)\times(W_0,\mbox{$\sum^{4n-1}$} \cup S^{4n-1}) \rightarrow Qp(k-1)\times(S^{4n-1}\times I,~S^{4n-1}\times\{0\}\cup S^{4n-1} \times\{1\}),
\]
with $I\times f|_{QP( k-1)\times S^{4n-1}}=$Identity.\\
The remainder of the proof is computing the obstruction $\sigma$ for this map to be a cobordism
and using this to determine if $Qp(k-1)\times\sum^{4n-1}$ is diffeomorphic to $Qp(k-1)\times S^{4n-1}$.\\\\
{\bf Claim:} $~~~~~~~~~~~~~$ $Ker( I\times f )_{\ast}= H_{\ast}( QP( k-1)) \times Ker(f_{\ast}).$\\\\
By the K$\ddot{\mbox{u}}$nneth formula, since $ H_{\ast}( QP( k-1))$ is torsion free then,
\[
H_{\ast}(QP( k-1) \times(W_0,\mbox{$\sum^{4n-1}$}\cup S^{4n-1}))\cong H_{\ast}( QP( k-1))\otimes H_{\ast} (W_0,\mbox{$\sum^{4n-1}$}\cup S^{4n-1})
\]
and $~~( I\times f )_{\ast}=I \otimes f_{\ast}$\\\\
therefore, $ Ker( I\times f )_{\ast}= H_{\ast}( QP( k-1)) \times Ker(f_{\ast}).$\\\\
Now, consider the commuative diagram induced by $f$:
\begin{eqnarray*}
 \longrightarrow H_{i}(\partial W_0)\longrightarrow~~~~~~~~~~~~~~~~H_{i}(W_0)\longrightarrow~~&H_{i}(W_0,\partial W_0)&\longrightarrow H_{i-1}(\partial W_0)\longrightarrow \hspace{3cm}\\
~\downarrow~~~~~~~~~~~~~~~~~~~~~~~~~~\downarrow~~~~~~~~~~~~~&f_{\ast}\downarrow&~~~~~~~~~~\downarrow\hspace{3cm}\\
\rightarrow H_{i}(\partial S^{4n-1}\times I)~~~~~~~\rightarrow ~H_{i}( S^{4n-1}\times I)\rightarrow ~&H_{i}( S^{4n-1}\times I,\partial S^{4n-1}\times I)&~\rightarrow H_{i-1}(\partial S^{4n-1}\times I)\rightarrow \hspace{3cm}
\\
\end{eqnarray*}
Notice that, $H_{i}(\partial S^{4n-1}\times I)\cong H_{i}(\partial W_0)~~\mbox{and}~~H_{i}( W_0) \cong 0 ~~\mbox{for}~~i\neq0,~2n $. We get
\begin{eqnarray*}
0\longrightarrow H_{2n}(W_0)\buildrel {\cong} \over\longrightarrow ~&H_{2n}(W_0,\partial W_0)&~\longrightarrow\hspace{3cm}\\
 ~&f_{\ast}\downarrow&~ \hspace{3cm}\\
~~0~~\longrightarrow ~~~0~~~ ~\longrightarrow~~&0&\longrightarrow~~~~0~~~\hspace{3cm}
\end{eqnarray*}
Hence $Ker(f_{\ast})\cong H_{2n}(W_0)$.\\\\
But $Ker({I\times f})_{\ast}=Ker({I\otimes f}_{\ast})= H_{\ast}(QP(k-1))\otimes Ker( f_{\ast}) \cong H_{\ast}(QP(k-1))\otimes H_{2n}(W_0),$ and
 $Ker({I\times f})_{\ast 2n+2k-2}=H_{2k-2}(QP(k-1))\otimes H_{2n}(W_0),$\\
 since $2k-2= 2(2s)-2\neq 0~\mbox{mod}~(4)$, $k$ is even, then $H_{2k-2}(QP(k-1))\cong 0$ and $Ker({I\times f})_{\ast 2n+2k-2}\cong 0$.\\
 Therefore, $\sigma(I\times f)=0$ and there exists an $h-$ cobrdism between $QP(k-1)\times \sum^{4n-1}$ and $QP(k-1)\times S^{4n-1}$, but
$k \geq 2$ and $n \geq 1$, then $(4k- 4)+(4n- 1)= 4(k+ n) - 5\geq 7.$
Hence, Smale's $h-$cobordism theorem can be applied and $QP(k-1)\times \sum^{4n-1}$ is
diffeomorphic to $QP(k-1)\times S^{4n-1}$. Applying Theorem 2.3, it follows that there is a
semi$-$free action of $S^3$ on some homotopy sphere $\sum^m$ with $\sum^{4n-1}$ as untwisted fixed
point set, where $m=4(n+k)-1.$\\
\begin{center}
{\Large\textbf{Open Problem :}}
\end{center}
Browder \cite{15} showed how to construct semi$-$free $S^1$ actions, with$\sum^r$ as
untwisted fixed point set, i.e., it�s normal bundle is trivial. He stated that he did not
know of any action with a twisted fixed point set. However Schultz \cite{11} used
complicated computations of homotopy groups, proved the following Theorem:\\
Let $k\geq 2$ be a positive integer. Then there exist infinitely many values of $n$ for which
$S^{2n}$ has a semi$-$free $S^1$ action with $S^{2(n-k+1)}$ as twisted fixed point set. In this work,
as in the work of Browder, we consider$\sum^r$ as untwisted fixed point set and $I$ raise
the following question: Does there exist smooth semi free $S^3$ actions on homotopy
spheres for which the fixed point set is twisted?\\\\
{\Large\textbf{Acknowledgement}}  $~$ I am very grateful to Professor Adil G. Naoum who
guided me in different topics of this paper.\\


\begin{thebibliography}{99}
\bibitem{1} A.Gleason, {\em Spaces with a compact Lie group of transformations}, proc.Amer.Math.Soc.I (1950), 35-43.
\bibitem{2} Eilenberg, Samuel; Steenrod, Norman , {\em Foundations of Algebraic Topology,} Princeton University Press, (1952), ISBN 978-0691079653.
\bibitem{3} F.Uchida, {\em Cobordism groups of semifree $S^1$ and $S^3$ actions}, Osaka.J.Math.7(1970),345-351.
\bibitem{4} H.T.Ku, {\em A note on Semi$-$free action on homotopy spheres}, Proc. Amer. Math. Soc.22(1969),
614-617.
\bibitem{5}H.T.Ku,and M.C.Ku,  {\em Free differentiable actions of $S^1$ and $S^3$ on homotopy spheres}. Proc. Amer .Math.Soc.25(1970),864-869.
\bibitem{66} I.H.Kaddoura, {\em De$-$suspension of free $S^3-$actions on Homotopy spheres,} reprint: cite as: http://arxiv.org/a/kaddoura$_{-}i_{-}1$.
\bibitem{6} J.Levine, {\em A classification of differentiable knots}, Ann.of Math.,82(1965), 15-50.
\bibitem{7} M.Kervaire, and J.Milnor, {\em Groups of homotopy spheres,} I.Ann. of Math. 77(1963), 504-537, 22.
\bibitem{8} M.J.Greenberg, {\em Lectures on algebraic topology}, W.A.Benjamin, Inc., New york, (1967).
\bibitem{9} P.E.Conner, and E.E.Floyd, {\em Differentiable periodic maps}; Springer Verlag, Acadamic press, Inc.Publishers, New york, (1964).
\bibitem{10} R. Brown, {\em Groupoids and Van Kampen�s theorem}, Proc. London Math. Soc.(3) 17(1967) 385-401.
\bibitem{11} R.Schultz, {\em Semifree actions with twisted fixed point sets}, Proc.Conf.on transformation groups, Springer Verlag, Berlin, Hidelberg, New York, (1972), 102-116.
\bibitem{12}S.P.Novikov,{\em Homotopically equivalent smooth manifolds}, Translations Amer.Math.Soc. 48(1965), 271-396.
\bibitem{13} S.Smale, {\em On the structure of manifolds,} Amer. J. Math. , 84(1962) pp.387�399
\bibitem{14} S.T.Hu, {\em Homology theory} , Holden-day, Inc.(1970).
\bibitem{15} W.Browder, {\em Surgery and the theory of differentiable transformation groups,} Proc. Conference on transformation groups (New Orleans,1967 ) New York,(1968), 1-46.
\bibitem{16} W.C.Hsiang, {\em A note on free differentiable actions of $S^1$ and $S^3$ on homotopy spheres}, Ann. of Math.83 (1966), 266-272.
\end{thebibliography}
\end{document}